\documentstyle{article}

\newtheorem{theorem}{Theorem}
\newtheorem{lemma}[theorem]{Lemma}
\newtheorem{conjecture}[theorem]{Conjecture}
\newtheorem{corollary}[theorem]{Corollary}
\newtheorem{fact}[theorem]{Fact}
\newtheorem{definition}{Definition}

\newtheorem{remark}{Remark}

\def\qedf{\hfill $\Box$}

\title{{\bf  On embedding well-separable graphs}}

\author{B\'ela Csaba
\thanks{Part of this research was done during the author's stay at Max-Planck-Institut f\"ur Informatik, 
Saarbr\"ucken, Germany} 
\thanks{Partially supported by the IST Programme of the EU under contract number IST-1999-14186 
(ALCOM-FT), and by OTKA T034475.}  
\\ Analysis and Stochastics Research Group of the \\ Hungarian Academy of Sciences\\ 
6720 Szeged, Aradi Vertanuk tere 1, Hungary \\ e-mail:bcsaba@math.u-szeged.hu}

\date{}

\begin{document}
\maketitle
\begin{abstract}
Call a simple graph $H$ of order $n$ {\it well-separable}, if by deleting a separator set
of size $o(n)$ the leftover will have components of size at most $o(n)$.
We prove, that  bounded degree well-separable spanning subgraphs are easy
to embed: for every $\gamma >0$ and positive integer $\Delta$ there exists an $n_0$ such 
that if $n>n_0$,  $\Delta(H) \le \Delta$ for a well-separable graph $H$ of order $n$ and  
$\delta(G) \ge (1-{1 \over 2(\chi(H)-1)} + \gamma)n$ for a simple
graph $G$ of order $n$, 
then  $H \subset G$. We extend our result to graphs with small band-width, too. 
\end{abstract}

\section{Notation} \label{jeloles}
In this paper we will consider only simple graphs.
We mostly use standard graph theory notation: we denote by $V(G)$ and $E(G)$ the vertex
and the edge set of the graph $G$, respectively. $deg_G(x)$ (or $deg(x)$) is the degree of  
the vertex $x \in V(G)$, $\delta(G)$ is the minimum
degree and  $\Delta(G)$ is the maximum degree. Denote  $deg_G(v,A)$ 
the number of  neighbors of $v$ in the set $A$.
We write $N_G(x)$ (or $N(x)$) for the neighborhood of the vertex $x \in V(G)$, hence, 
$deg_G(x)=|N_G(x)|$.
$N_G(U)=\cup_{x \in U}N(x)$ for a set $U \subset V(G)$. $N_G(v,A)$  is the set of
neighbors of $v$ in $A$. Set $e(G)=|E(G)|$ and $v(G)=|V(G)|$. If $A$ and $B$ are disjoint
subsets of $V(G)$, then we denote by $e(A,B)$ the number of edges with one endpoint in $A$
and the other in $B$. We write $\chi(G)$ for the chromatic number of $G$. If $A$ is a subset
of the vertices of $G$, we write $G-A$ for the graph induced by the vertices of $V(G)-A$.

If $G$ has a subgraph isomorphic to $H$, then we write $H \subset G$. In this case
we sometimes call $G$ the {\it host graph}. We say that $G$ has an
$H$--factor if there are $\lfloor v(G)/v(H) \rfloor$ vertex-disjoint copies of $H$ in $G$ 
(this notion is somewhat different from the common one: we don't need
that $v(G)$ is a multiple of $v(H)$).
Throughout the paper we will apply the relation 
``$\ll$'': $a \ll b$ if $a$ is sufficiently smaller than $b$.

\section{Introduction} \label{bevezetes}

In this paper we consider a problem in extremal graph theory.
Before getting on the subject of our result let us take a short
historical tour in the field.

One of the main results of the area is Tur\'an's Theorem:

\begin{theorem}[Tur\'an 1941~\cite{T}] \label{Turan}
If $G$ is a graph on $n$ vertices, and
$$e(G)>\left(1-{1 \over r-1}\right){n^2 \over 2},$$
then $K_r \subset G$.
\end{theorem}

Another milestone in extremal graph theory is the following theorem:

\begin{theorem}[Erd\H{o}s--Stone--Simonovits 1946/1966~\cite{ESt,ESi}] \label{E-S-S}
For every graph $H$ and every real $\varepsilon >0$ there exists an $N=N(H,\varepsilon)$
such that if $G$ is a graph on $n>N$ vertices, and
$$e(G)>\left(1-{1 \over \chi(H)-1}+\varepsilon \right){n^2 \over 2},$$
then $H \subset G$.
\end{theorem}

The deep result of Hajnal and Szemer\'edi shows that when we are looking
for a $K_r$--factor in a graph, the situation is different.

\begin{theorem}[Hajnal--Szemer\'edi 1969~\cite{HSz}] \label{H-Sz}
If $G$ is a graph of order $n$ and $\delta(G) \ge (1-1/r)n$, then $G$ has
a $K_r$--factor.
\end{theorem}

There are two important changes in the formulation of the  above result: first, it is not 
sufficient
to bound the number of edges anymore -- we need a lower bound on the minimum degree
of the host graph. Second, that $1/(r-1)$ changed to $1/r$.

The following results were conjectured by Alon and Yuster~\cite{AY92, AY96}, 
and proved by Koml\'os, S\'ark\"ozy and Szemer\'edi:

\begin{theorem}[Koml\'os--S\'ark\"ozy--Szemer\'edi 2001~\cite{KSSz}]\label{K-S-Sz}
{\bf Part 1:} For every graph $H$ there is a constant $K$ such that if $G$ is a graph on $n$ vertices, then  
$$\delta(G) > \left(1-{1 \over \chi(H)} \right)n$$
implies that there is a union of vertex disjoint copies of $H$ covering all
but at most $K$ vertices of $G$.

\smallskip

\noindent {\bf Part 2:} For every graph $H$ there is a constant $K$ such that if $G$ is a graph on $n$ 
vertices, then  $$\delta(G) > \left(1-{1 \over \chi(H)} \right)n + K$$
implies that $G$ has an $H$--factor.
\end{theorem}

These theorems show that the chromatic number is a crucial parameter in classical extremal
graph theory. However, it is easy to come up
with examples when the maximum degree turns out to be much more important.
We give one possible set of examples for this fact.
Let $\{H_d\}_{d>2}$ be a family of random bipartite graphs with equal color
classes of size $n/2$ that are obtained as the union of $d$ random 1--factors.
Let $r$ be an odd positive integer, and consider the graph $G$ of order $n$ 
having $r$  independent sets of equal size, and all the edges between any 
two independent sets.
By a standard application of the probabilistic method one can prove that for a given $r$
if $d$ is large enough ($d={\rm constant} \cdot r$ is sufficient), then $H_d \not\subset G$. 
Since $H_d$ is bipartite
for every $d$, this proves, that the critical parameter for embedding expanders
cannot be the chromatic number. (Although, the chromatic number still has
a role, see~\cite{Cs03}.)
One may think, that the main reason of this fact is that $H_d$ is an expander
graph with large expansion rate. 

We show, that if a graph is "far from being an expander", then again, the chromatic
number comes into picture. First, let us define what we mean on "non--expander" graphs.

\begin{definition}\label{well-separable}
Let $H$ be a graph of order $n$. We call $H$ well-separable if there is a subset
$S \subset V(H)$ of size $o(n)$ such that all components of $H-S$ are of size $o(n)$.
\end{definition}

We call $S$ the separator set, and write $C_1,C_2,\ldots,C_t$ for the components
of $H-S$.
Note, that if $H$ is an expander graph, then it is not well--separable.
We will show the following property of well--separable graphs.

\begin{theorem}\label{tetel}
For every $\gamma >0$, positive integers $\Delta$ and $k$ there exists an $n_0$ such 
that if $n>n_0$, $\chi(H)\le k$, $\Delta(H) \le \Delta$ for a well-separable graph $H$ 
of order $n$ and  $\delta(G) \ge (1-{1 \over 2(k-1)} + \gamma)n$ for a simple
graph $G$ of order $n$, then  $H \subset G$.
\end{theorem}

Observe, that trees are well--separable graphs. A conjecture of Bollob\'as~\cite{B} (proved by
Koml\'os, S\'ark\"ozy and Szemer\'edi~\cite{KSSz95}) states that trees of bounded degree can be 
embedded into graphs of minimum degree $(1/2 +\gamma)n$ for $\gamma >0$. Since every tree is bipartite, this
result is a special case of Theorem~\ref{tetel}. (Recently Koml\'os, S\'ark\"ozy and Szemer\'edi extended 
their result
for trees of maximum degree as large as $c {n \over \log{n}}$~\cite{KSSz2}.)

Our proof of Theorem~\ref{tetel} uses the Regularity Lemma of Szemer\'edi~\cite{Sz76}
(sometimes called Uniformity Lemma). In the next section we will give a brief survey on 
this powerful tool, and related results. For more information see e.g.,~\cite{KS93,K}. 
We will prove Theorem~\ref{tetel} in the fourth section, and then prove a strengthened version of it, too. 
In the fifth section we will investigate
the case of graphs with small band-width.

\section{A review of tools for the proof}

We introduce some more notation first. The {\it density}  between disjoint sets $X$ and $Y$ 
is defined as:
$$d(X,Y)={{e(X,Y)}\over {|X||Y|}}.$$

We need  the following
definition to state the Regularity Lemma.

\begin{definition}[Regularity condition] Let $\varepsilon >0$. A pair
  $(A,B)$      of     disjoint      vertex sets     in      $G$     is
  \mbox{$\varepsilon$-regular}  if  for every  $X  \subset  A$ and  $Y
  \subset B$, satisfying
$$|X|>\varepsilon |A|,\ |Y|>\varepsilon |B|$$
we have
$$|d(X,Y)-d(A,B)|<\varepsilon.$$
\end{definition}

We will employ the fact that if $(A,B)$ is an $\varepsilon$--regular pair as above, and
we place $constant \cdot \varepsilon |A|$ new vertices into $A$, the resulting
pair will remain $\varepsilon'$-regular, with a somewhat larger $\varepsilon'$ than $\varepsilon$, depending
on the constant.

An important property of regular pairs is the following:

\begin{fact}\label{metszes}
Let $(A,B)$ be an $\varepsilon$--regular pair with density $d$. Then for any $Y \subset B$,
$|Y|>\varepsilon |B|$, we have
$$ |\{x \in A:deg(x,Y)\le (d-\varepsilon)|Y|\}|\le \varepsilon |A|.$$
\end{fact}

We will use the following form of the Regularity Lemma:

\begin{lemma}[Degree Form] \label{rl} For every $\varepsilon>0$ there is an $M=M(\varepsilon)$ 
such
that if $G=(V,E)$ is any graph and $d\in [0,1]$ is any real number, then there is a partition 
of the
vertex set $V$ into $\ell+1$ clusters $V_0, V_1,\ldots,V_\ell$, and there is a subgraph $G'$ 
of $G$
with the following properties:  
\begin{itemize} 

\item $\ell \le M$, 

\item $|V_0|\le \varepsilon |V|$,

\item all clusters $V_i$, $i\ge 1$, are of the same size $m$ $
($and therefore $m \le \lfloor {|V|\over \ell}\rfloor<\varepsilon |V|)$, 

\item $deg_{G'}(v)>deg_{G}(v)-(d+\varepsilon)|V|$ for all $v\in V$, 

\item $V_i$ is an independent set in $G'$ for all $i\ge 1$, 

\item all pairs $(V_i,V_j)$, $1\le i <j \le \ell$, are $\varepsilon$-regular, each with 
density either 0 or at least $d$ in $G'$.  

\end{itemize} 

\end{lemma}

\noindent Often  we  call  $V_0$  the  {\it exceptional  cluster}.
In the rest of the paper we assume that $0< \varepsilon \ll  d \ll  1$.

\begin{definition}[Reduced graph]
Apply Lemma~\ref{rl} to the graph $G=(V,E)$ with parameters $\varepsilon$ and $d$, 
and denote the clusters of the resulting partition by $V_0, V_1,\ldots,V_\ell$, $V_0$ being
the exceptional cluster. We construct a new graph $G_r$, the reduced graph of $G'$ in
the following way: The non-exceptional clusters of $G'$
are the vertices of the reduced graph (hence $|V(G_r)|=\ell$). We connect two vertices of
$G_r$ by an edge if the corresponding two clusters form an $\varepsilon$-regular
pair with density at least $d$.
\end{definition}

The following corollary is immediate:

\begin{corollary}\label{redukaltfok}
Let $G=(V,E)$ be a graph of order $n$ and  $\delta(G) \ge  cn$ for some $c>0$, and let $G_r$ 
be the reduced graph of $G'$
after applying Lemma~\ref{rl}  with parameters $\varepsilon$ and $d$. 
Then $\delta(G_r) \ge (c - \theta)\ell$,
where $\theta=2\varepsilon+d$.
\end{corollary} 

A  stronger one-sided property of regular pairs is super-regularity:

\begin{definition}[Super-Regularity condition] Given a graph $G$ and
  two disjoint subsets   $A$ and  $B$ of its vertices, the pair $(A,B)$
  is   \mbox{$(\varepsilon,\delta)${\rm    -super-regular}},   if   it   is
  $\varepsilon$-regular and furthermore,
$$deg(a)>\delta |B|,{\rm \ for\  all}\  a\in A,$$
and
$$deg(b)>\delta |A|,{\rm \ for\  all}\  b\in B.$$
\end{definition}

Finally, we formulate another important tool of the area:

\begin{theorem}[Blow-up Lemma~\cite{KSSz97,KSSz98}]\label{blow-up}
Given a graph $R$ of order $r$ and positive parameters $\delta, \Delta$, there exists
a positive $\varepsilon=\varepsilon(\delta, \Delta, r)$ such that the following holds:
Let $n_1,n_2, \ldots, n_r$ be arbitrary positive integers and let us replace the vertices
$v_1, v_2, \ldots, v_r$ of $R$ with pairwise disjoint sets $V_1,V_2, \ldots, V_r$ of
sizes $n_1,n_2, \ldots, n_r$ (blowing up). We construct two graphs on the same vertex
set $V=\cup V_i$. The first graph $F$ is obtained by replacing each edge $\{v_i,v_j\}$
of $R$ with the complete bipartite graph between $V_i$ and $V_j$. A sparser graph $G$
is constructed by replacing each edge $\{v_i,v_j\}$ arbitrarily with an 
$(\varepsilon,\delta)$--super--regular pair between $V_i$ and $V_j$. If a graph $H$ with 
$\Delta(H) \le \Delta$
is embeddable into $F$ then it is already embeddable into $G$. 
\end{theorem}
 
\begin{remark}[Strengthening the Blow-up Lemma~\cite{KSSz97}]\label{erosites}
Assume that $n_i \le 2 n_j$ for every $1 \le i,j \le r$. Then we can strengthen the lemma:
Given $c>0$ there are positive numbers 
$\varepsilon=\varepsilon(\delta,\Delta,r,c)$ and $\alpha=\alpha(\delta,\Delta,r,c)$ such that
the Blow-up Lemma remains true if for every $i$ there are certain vertices $x$ to be embedded
into $V_i$ whose images are a priori restricted to certain sets $T_x \subset V_i$ provided that
\begin{itemize}
\item[](i) each $T_x$ within a $V_i$ is of size at least $c|V_i|,$
\item[](ii) the number of such restrictions within a $V_i$ is not more than $\alpha|V_i|.$
\end{itemize}
\end{remark}

\section{Proof of Theorem~\ref{tetel}}

The proof goes along the following lines:

\begin{itemize}

\item[](1) Find a special structure in $G$ by the help of the Regularity Lemma and the
Hajnal--Szemer\'edi Theorem (Theorem~\ref{H-Sz}).

\item[](2) Map the vertices of $H$ to clusters of $G$ in such a way that if $\{x,y\}\in E(H)$,
then $x$ and $y$ are mapped to neighboring clusters; moreover, these clusters will form
an $(\varepsilon, \delta)$--super--regular pair for all, but at most $o(n)$ edges.

\item[](3) Finish the embedding by the help of the Blow-up Lemma.

\end{itemize}

\subsection{Decomposition of $G$} \label{struktura}

In this subsection we will find a useful decomposition of $G$.

First, we apply the Degree Form of the Regularity Lemma with parameters $\varepsilon$ and
$d$, where $0<\varepsilon \ll d \ll \gamma < 1$.
As a result, we  have $\ell+1$ clusters, $V_0, V_1, \ldots, V_{\ell}$, where $V_0$ is
the exceptional cluster of size at most $\varepsilon n$, and all the others have the same 
size $m$. We deleted only a small number of edges, and now all the $(V_i,V_j)$ pairs
are $\varepsilon$--regular, with density 0 or larger than $d$. By Corollary~\ref{redukaltfok}
we will have that $\delta(G_r) \ge (1- {1 \over 2(k-1)}+\gamma')\ell$, where
$\gamma'=\gamma -d-2\varepsilon >0$.

Applying Theorem~\ref{H-Sz}, we have a $K_{k}$--factor in $G_r$. It is possible,
that at most $k-1$ clusters are left out from this $K_k$--factor -- such clusters
are put into $V_0$.
It is easy to transform the $\varepsilon$--regular pairs inside this $K_k$--factor into
super--regular pairs: given a $\delta$ with $\varepsilon \ll \delta \ll d$
we have to discard at most $\varepsilon m$ vertices from a cluster to
make a regular pair $(\varepsilon,\delta)$--super--regular. In a $k$--clique a cluster has $k-1$
other adjacent clusters in $G_r$. Hence, it is enough to discard at most
$(k-1)\varepsilon m$ vertices from every cluster, and arrive to the
desired result. Note, that now the pairs are $\varepsilon'$--regular, with
$\varepsilon'< 2 \varepsilon$; for simplicity, we will use the letter $\varepsilon$
in the rest of the paper. We will discard the same number of vertices
from every non--exceptional cluster, and get, that all the edges of $G_r$
inside the cliques of the $K_k$--factor are $(\varepsilon,\delta)$--super-regular pairs.
For simplicity we will still denote the common cluster size by $m$ in $G_r$.
The discarded vertices are placed into $V_0$; now $|V_0| \le (2k-1) \varepsilon n$.

Our next goal is to distribute the vertices of $V_0$ among the non--exceptional
clusters so as to preserve super--regularity within the cliques of the $K_k$--factor. We also 
require
that the resulting clusters should have about the same size.

For a cluster $V_i$ in $G_r$ denote $clq(V_i)$ the set of the clusters of $V_i$'s
clique in the $K_k$--factor, but without $V_i$ itself. Hence, $V_i \not\in clq(V_i)$, 
and $|clq(V_i)|=k-1$ for every
$V_i \in V(G_r)$.

Recall, that every cluster in $G_r$ has the same size, $m$.
We want to distribute the vertices of $V_0$ evenly among the clusters of $G_r$: 
we will achieve that $||V_i|-|V_j||<4k \varepsilon m$ for every $1 \le i,j \le \ell$
after placing the vertices of $V_0$ to non--exceptional clusters.
Besides, we  require that if we put a vertex $v \in V_0$ into $V_i \in V(G_r)$, then 
$deg(v, V_j) \ge \delta m$ for every $V_j \in clq(V_i)$.  

So as to satisfy the above requirement, let us define an auxiliary bipartite graph \\
$F_1=F_1(V_0,V(G_r),E(F_1))$. That is, the color classes of $F_1$ are $V_0$ and the
set of the non--exceptional clusters. We draw a $\{v, V_i\}$ edge for $v \in V_0$ and
$V_i \in V(G_r)$ if $deg_G(v, V_j) \ge \delta m$ for every $V_j \in clq(V_i)$.

Set $\gamma''=k(\gamma-2(\varepsilon+d))$.
The following lemma is crucial in distributing $V_0$.

\begin{lemma} \label{LP}
$deg_{F_1}(v) \ge (1/2 +\gamma'')\ell$ for every $v \in V_0$.
\end{lemma}

\noindent {\bf Proof:}
Consider an arbitrary $v \in V_0$. Then we can partition the set of $k$--cliques of the
$K_k$--factor into
$k+1$ pairwise disjoint sets $A_0, A_1, \ldots, A_k$. A clique $Q$ is in $A_j$ if $v$
has at least $\delta m$ neighbors in exactly $j$ clusters of $Q$. Set $a_j=k|A_j|/\ell$ for
every $0 \le j \le k$, that is, $a_j$ is the proportion of cliques in $A_j$. Clearly,
$\sum_j a_j=1$. There are at most $\delta n$ edges connecting  $v$ to clusters not adjacent to $v$ in $F_1$. 
Hence, by the minimum degree condition, 
$1/k \sum_j {j a_j} \ge \delta(G_r)/ \ell -\delta$. Notice, that if there are at most $k-2$ clusters in 
a clique in which $v$ has at least  
$\delta m$ neighbors, then $v$ is not adjacent to any clusters  of that clique
in $F_1$. There are two possibilities left: $v$ has one neighbor in a clique in $F_1$,
or it is connected to all the clusters in $F_1$, depending on whether it has large 
enough degree
to $k-1$ or $k$ clusters of that clique. 
Putting these together, the solution of the following linear program is a lower
bound for $deg_{F_1}(v)/ \ell$:

\begin{eqnarray}
\nonumber &{\displaystyle \sum_{j=0}^k{a_j}=1} \ \ {\rm and} \ \ 
{\displaystyle \sum_{j=0}^k{j a_j} -z}=k({2k-3 \over 2k-2}+\gamma-2(\varepsilon+d)) \\
\nonumber &{\rm where} \ \ a_j, z \ge 0 \\
\nonumber &\min \{ {a_{k-1} \over k} +a_k \}
\end{eqnarray}

Let $A$ be the coefficient matrix of the two equalities above, i.e.,  
$$A= \pmatrix{1 & 1 & 1 &  \ldots & 1 & 1 & 0 \cr 0 & 1 & 2 & \ldots & k-1 & k & -1}.$$
Let $a^T=(a_0,a_1, \ldots, a_k,z)$, $b^T=(1, k(2k-3)/(2k-2)+\gamma'')$, and
$c^T=(0,0,\ldots,0,1/k,1,0).$
Then the dual of the linear program above is:

\begin{eqnarray}
\nonumber &A^T u \le c \\ 
\nonumber &\max \{b^T u \} 
\end{eqnarray}

It is easy to check that $u_1=2-k$ and $u_2 = {k-1 \over k}$ is a feasible
solution (in fact the optimal solution as well), and therefore 
$\max{b^T u} \ge 1/2 + \gamma''$. \hfill \qedf

\bigskip

Applying the lemma above it is easy to distribute the vertices of $V_0$ 
evenly, without violating our requirement. For every $v \in V_0$ randomly
choose a neighboring cluster in $F_1$, and put $v$ into that cluster.
Since $deg_{F_1}(v) \ge (1/2 + \gamma'')\ell$, with very high probability
(use eg., Chernoff's bound) no cluster will get more than $2 |V_0|/\ell$
new vertices from $V_0$. Hence, we have that  $||V_i|-|V_j||<4k \varepsilon m$ 
for every $1 \le i,j \le \ell$.

\subsection{Assigning the vertices of $H$} \label{lekepezo}

In this subsection we will map the vertices of $H$ to clusters of $G_r$.
We will heavily use the fact that $H$ is $k$--colorable.

Fix an arbitrary $k$--coloration of $H$. For an arbitrary set $A$, denote
$A^1,A^2,\ldots,A^k$ the color classes determined by this $k$--coloration.

Recall, that $S$ is the separator set of $H$ and $C_1, C_2, \ldots, C_t$ are the components of $H-S$.
We will map $S$ and $C_1,C_2,\ldots, C_t$ by the randomized procedure below. \\

\noindent {\bf Mapping algorithm} 

\medskip

Input: the set $A$

\begin{itemize}

\item Pick a clique $Q=\{Q_1,Q_2,\ldots,Q_k\}$ in the cover of $G_r$ randomly, uniformly.

\item Pick a permutation $\pi$ on $\{1,2,\ldots,k\}$ uniformly at random.

\item Assign the vertices of $A^i$ to the cluster $Q_{\pi(i)}$ for every $1 \le i \le k$.

\end{itemize}

Repeating this algorithm for $S$ and all the components in $H-S$, we will have, that
the number of vertices of $H$ assigned to a cluster are almost the same: with probability
tending to 1, the  difference between the number of assigned vertices to a cluster and the 
cluster size $m$ will be at most $o(n)$. This follows easily from a standard application of 
Chebyschev's inequality:

\begin{lemma}\label{Csebisev}
With positive probability the mapping algorithm assigns $n/\ell \pm \varepsilon m/\ell$ vertices of $H$ to every 
cluster of $G_r$.
\end{lemma}

\noindent {\bf Proof:}
Let $V_t$ be an arbitrary cluster of $G_r$. The above mapping algorithm is a randomized procedure, hence, 
the number 
of vertices of $H$ assigned to $V_t$ is a random variable. Denote this random variable by $Z$. Let us 
define $n$ indicator
random variables $ \{Z_i \}_1^n$, where $Z_i=1$ if and only if $x_i$ (the $i$th vertex of $H$) 
is assigned to $V_t$ by
the mapping algorithm. Notice, that these indicator variables follow the same distribution. Clearly, 

$$Z=\sum_{i=1}^n{Z_i},$$

hence, ${\rm E}(Z)=n / \ell$ and 

$${\rm Var}(Z)=\sum_{i=1}^n{\rm Var}(Z_i)+ \sum_{i \neq j}^n{({\rm E}(Z_i Z_j)-
{\rm E}(Z_i){\rm E}(Z_j))}.$$ 

If $Z_i$ and $Z_j$ are independent, then ${\rm E}(Z_i Z_j)={\rm E}(Z_i){\rm E}(Z_j)$. 
Since $Z_i$ and
$Z_j$ are independent if they belong to different components, we can give a trivial upper bound on 
${\rm Var}(Z)$:

$${\rm Var}(Z) \le n{\rm Var}(Z_1) +o(n) n.$$

Let us apply Chebyschev's inequality for $Z$: $${\rm Pr}(|Z -n/ \ell| \ge \lambda {\rm D}(Z)) \le 
{1 \over \lambda^2}$$
with ${\rm D}(Z) = \sqrt{{\rm Var}(Z)}$.
Observe, that if $1 /\lambda^2 < 1 /\ell$ then we can guarantee that the mapping algorithm assigns 
$n/\ell \pm \lambda {\rm D}(Z)$ 
vertices of $H$ to every cluster with positive probability. We set $\lambda=\sqrt{2 \ell}$. Since  
${\rm D}(Z) = o(n)$, we have that 
$\lambda {\rm D}(Z) < \varepsilon m/\ell$ if $n$ is large enough. 
\hfill \qedf

\medskip

Recall, that for applying the Blow--up Lemma, it is necessary to map adjacent
vertices in $H$ to adjacent clusters in $G_r$.
For $x \in V(H)$ let $\kappa(x)$ denote the cluster to which $x$ is assigned.
After randomly assigning $S$ and $C_1,C_2, \ldots, C_t$, we have that if $\{x,y\} \in H$
and $x,y \in S$ or $x,y \in C_j$ for some $1 \le j \le t$, then 
$\{ \kappa(x),\kappa(y) \} \in E(G_r)$. On the other hand, there is no guarantee that 
a vertex in $S$ and a vertex in some component of $H-S$ are 
assigned to adjacent clusters, even when they are adjacent in $H$. 

Therefore, we have to reassign a small subset of $V(H)$. We will see that 
no vertex which is at distance larger than $k$ from $S$ will change its place,
and vertices of $S$ will not be reassigned. 
Consider an arbitrary component $C_j$. Set $B=N(S) \cap C_j$, and $B_p=B \cap C_j^p$ 
for every $1 \le p \le k$.  By the algorithm below we will define $B'_p$, the subset 
of $C_j^p$ which will be reassigned. 

\begin{itemize}

\item[]{\it Step 1.} Set $B'_k=B_k$, and $i=1$

\item[]{\it Step 2.} Set $B'_{k-i}=B_{k-i} \cup \bigcup_{p=0}^{i-1}
(N(B'_{k-p})\cap C_j^{k-i})$

\item[]{\it Step 3.} If $i < k-1$, then set $i \leftarrow i+1$, and go back to {\it Step 2.}

\end{itemize}

Informally, when we determine which vertices to reassign from $C_{k-i}$, we take
into account all the neighbors of $B'_p$ with $p>k-i$, and $B_{k-i}$ itself.
It is important, that we proceed backwards, that is, we specify the vertices to be reassigned 
starting from the last, the $k$th color class. 
Note, that the vertices of $\cup_{p=1}^{k}B'_p$ are at distance at most $k$ from $S$.
Hence, $|\cup_{p=1}^{k} B'_p|< \Delta^k |S|=o(n)$.

Now we have the sets $\{B'_p\}$. First we will find a new cluster for $B'_1$:
Take an arbitrary cluster $W_1$ from the set 

$$\bigcap_{p=2}^k N(\kappa(S^p)) \cap \bigcap_{p=2}^k N(\kappa(B'_p)),$$

and assign the vertices of $B'_1$ to the cluster $W_1$.

Then we choose $W_2$ for $B'_2$ from the set 

$$\bigcap_{p \not=2} N(\kappa(S^p)) \cap \bigcap_{p=3}^k N(\kappa(B'_p)) \cap N(W_1),$$

and assign the vertices of $B'_2$ to the cluster $W_2$.

In general, assume that we have the clusters $W_1, W_2, \ldots, W_{i-1}$ for some $i \le k$. 
Then we choose 
$W_i$ for  $B'_i$ from the set

$$\bigcap_{p \not=i} N(\kappa(S^p)) \cap \bigcap_{p=i+1}^k N(\kappa(B'_p)) \cap 
\bigcap_{p=1}^{i-1} N(W_p),$$

and assign the vertices of $B'_i$ to the cluster $W_i$.

Observe, that this way $W_i$ ($1 \le i \le k$) is chosen from a non--empty set, since it 
comes from the common neighborhood of $2k-2$ clusters, and this neighborhood is of size 
at least $\gamma' \ell$ by the minimum degree condition of $G$. 

By the help of the above reassigning procedure we achieved, that  adjacent vertices of $H$ are assigned to adjacent
clusters of $G_r$. Let us denote the set of vertices of $H$ assigned to cluster $V_i$ by $L_i$ for 
every $1 \le i \le \ell$.
Our next goal is to make $|L_i|=|V_i|$.

\subsection{Achieving $|V_i|=|L_i|$} \label{kiegyenlito}

We have, that if $\{ x,y \} \in E(H)$, then $\{ \kappa(x), \kappa(y) \} \in E(G_r)$.
Moreover, the $\{ \kappa(x), \kappa(y) \}$ edges are super--regular pairs for all,
but at most $o(n)$ edges in $E(H)$. 

Still, we cannot apply the Blow--up Lemma, since $|V_i|=|L_i|$ is not necessarily 
true for every $1 \le i \le \ell$. What we know for sure is that 
$||V_i|-|L_i||<5 \varepsilon k m$, because these differences were at most $o(n)$
after the random mapping algorithm of the previous subsection, and distributing the vertices of $V_0$ had 
contribution at most
$4 k \varepsilon m$ for every $1 \le i \le \ell$ (we refer to Subsection~\ref{struktura}), and we relocated 
$o(n)$ vertices
in the previous subsection.

We will partition the clusters of $G_r$ into three disjoint sets: $V_<, V_=$ and
$V_>$. If $|V_i| < |L_i|$, then $V_i \in V_<$; if $|V_j|=|L_j|$, then $V_j \in V_=$,
and we put $V_p$ into $V_>$ if $|V_p|>|L_p|$. Clearly, it is enough to replace
at most $5k \varepsilon n$ vertices of $G$ so as to achieve $|V_i|=|L_i|$ for every
$1 \le i \le \ell$, while preserving regularity for the edges of $G_r$. But we
need super--regular pairs for the edges of the $k$--cliques of the $K_k$--factor, hence, 
a straightforward relocation of some vertices of $G$ is not helpful. Instead, we
will apply an idea similar to what we used for distributing the vertices of $V_0$.

First, we define a directed graph $F_2$: the vertices of $F_2$ are the clusters
of $G_r$, and $(V_i,V_j) \in E(F_2)$ if $(V_i, V_p) \in E(G_r)$ for every 
$V_p \in clq(V_j)$. 
We will have that the out--degree of every cluster is at least $(1/2 +\gamma'')\ell$
by considering the  linear program of Subsection~\ref{struktura}. 
Since $\delta(G_r) \ge ({2k-3 \over 2k-2}+\gamma')\ell$, it is easy to see that any
$k-1$ clusters have at least $(1/2 +\gamma')\ell$ common neighbors. That is, the in--degree
of $F_2$ is at least $(1/2 +\gamma')\ell$. Therefore, there is a large number  
-- at least $(\gamma'+\gamma'')\ell$ -- of 
directed paths of length at most two between any two clusters in $F_2$.

Let $V_i \in V_<$ and $V_j \in V_>$ be arbitrary clusters. If $(V_j, V_i) \in E(F_2)$, then
we can directly place a vertex from $V_j$ into $V_i$ which has at least $\delta m$
neighbors in $V_s$ for every $V_s \in clq(V_i)$ (and most of the vertices have actually at least $d m$ neighbors,
since $d$ is the lower bound for the density of regular pairs). If there is no such edge,
then there are several different directed paths of length
two from $V_j$ to $V_i$. These paths differ in their "center" cluster. Assume that $V_p$
is such a cluster, i.e., $(V_j,V_p)$ and $(V_p,V_i)$ are edges in $F_2$. It is useful
to choose $V_p$ randomly, uniformly among the possible "center" clusters. 

Take any vertex
$v \in V_j$ which has  at least $\delta m$ neighbors in $V_s$ for every $V_s \in clq(V_p)$, and put it into $V_p$.
Then choose any vertex from $V_p$ which has at least $\delta m$ neighbors in $V_t$ for every $V_t \in clq(V_i)$, 
and put it into $V_i$. As a result, we decreased $||V_j|-|L_j||$ and $||V_i|-|L_i||$, while
$||V_p|-|L_p||$ did not change. Now, by the remark after the definition of a regular pair
it is clear that if we make all $|V_i|=|L_i|$ this way, we will preserve regularity
and super--regularity as well.

\subsection{Finishing the proof}

Now we are prepared to prove Theorem~\ref{tetel}.

We have to check if the conditions of the Blow--up Lemma are satisfied.
There are $o(n)$ edges of $E(H)$ which are  problematic: those edges
having their endpoints in clusters which do not constitute a super--regular pair.
Denote the set of these edges by $E'$. Suppose that $x$ is
a vertex which occurs in some edges of $E'$. It can have neighbors assigned to
at most $2k-2$ clusters $V_{x_1},V_{x_2},\ldots,V_{x_{2k-2}}$. Since $(\kappa(x),V_{x_i})$
is a regular pair for every $1 \le i \le 2k-2$, there is a set $T_x \subset \kappa(x)$ of 
size at least $(1-(2k-2)\varepsilon)m$ (by Fact~\ref{metszes} and applying induction), all the 
vertices of which have at least
$(d-\varepsilon)^{2k-2}m > \delta m$ neighbors in $V_{x_i}$ for every $1 \le i \le 2k-2$. 
$T_x$ will be the set 
to which $x$ is restricted. Since $|E'|=o(n)$, the number of restricted vertices is small
enough, and therefore we can apply the strengthened version of the Blow-up Lemma.
\hfill \qedf

\subsection{Strengthening Theorem~\ref{tetel}}

We begin this subsection with a definition.

\begin{definition} Let $0< \alpha <1$. We call a graph $H$ on $n$ vertices $\alpha$--separable, if there is a set 
$S \subset V(H)$ of size
at most $\alpha n$ such that all components of $H-S$ are of size at most $\alpha n$. 
\end{definition}

Obviously, given some $0< \alpha <1$ if $H$ is well--separable and $|V(H)|$ is large enough, then $H$ is 
$\alpha$--separable as well.
On the other hand, if $\alpha$ is small enough, then we can substitute well--separability by $\alpha$--separability:

\begin{theorem}\label{tetel2}
For every $\gamma >0$, positive integers $\Delta$ and $k$ there exists an $n_0$ and an $\alpha$ such 
that if $n>n_0$, $\chi(H)\le k$, $\Delta(H) \le \Delta$ for an $\alpha$--separable graph $H$ 
of order $n$ and  $\delta(G) \ge (1-{1 \over 2(k-1)} + \gamma)n$ for a simple
graph $G$ of order $n$, then  $H \subset G$.
\end{theorem}
  
\noindent {\bf Proof (sketch):} We will apply the same method for embedding $\alpha$--separable graphs. First, we 
decompose $G$ by the help of the
Regularity Lemma and the Hajnal--Szemer\'edi Theorem. Then distribute the vertices of $H$ among the clusters of 
$G_r$, finally, apply the Blow-up Lemma
for finishing the embedding. Since $S$ and the components of $H-S$ can be much larger now, we have to be  
careful at certain points. We will
pay attention only to these points.

Given $\gamma, \Delta$ and $k$, we can determine $\alpha$: Proving Lemma~\ref{Csebisev} for $\alpha$--separable 
graphs 
we will have that ${\rm Var}(Z) \le n {\rm Var}(Z_1) + \alpha n^2$, hence, ${\rm D}(Z) \le 
\sqrt{2 \alpha}n$. Set $\lambda= \sqrt{2 \ell}$, and 
choose $\alpha$ so that
$${\rm Pr}(|Z -n/ \ell| \ge \varepsilon m) \le {1 \over \lambda^2}= {1 \over 2 \ell}. \quad\quad\quad \ (1)$$
It is easy to check that if $\alpha \le \varepsilon^2/(4 \ell^3)$ then $\lambda {\rm D}(Z) \le \varepsilon m$ 
and inequality (1) is satisfied. 

After the random mapping algorithm we have to reassign some vertices so as to get that adjacent vertices of $H$ are 
assigned to adjacent clusters
of $G_r$. At this point we may reassign as many as $\Delta^k |S| \le \Delta^k \alpha n$ vertices of $H$. Our second 
criteria for $\alpha$ is that  
$\Delta^k \alpha n$ should be less than $\varepsilon m$. Other parts of the proof work smoothly not just for 
well--separable but for $\alpha$--separable graphs as well.

Therefore, if 
$$\alpha \le \max \{{\varepsilon^2 \over 4 \ell^3}, {\varepsilon \over \ell \Delta^k}\}=
{\varepsilon^2 \over 4 \ell^3},$$ then we can embed $H$ into $G$.
\hfill \qedf

\section{On graphs with small band-width}

Another notion, which measures the "non-expansion" of  graphs is band-width. Let us denote the band-width of 
a graph  $G$ by $bw(G)$. 
Notice, that there are well--separable graphs with large band-width: consider $K_{1,n-1}$, the star on $n$ 
vertices.  Obviously, it is a well--separable 
graph, on the other hand its band-width is $n/2$.

The following is conjectured by Bollob\'as and Koml\'os (see e.g., in~\cite{KS93}):

\begin{conjecture}[Bollob\'as-Koml\'os] \label{BK}
For every $\gamma >0$ and positive integers $r$ and $\Delta$, there is a $\beta>0$ and an $n_0$ such that if 
$|V(H)|=|V(G)|=n \ge n_0$,
$\chi(H) \le r$, $\Delta(H) \le \Delta$, $bw(H)< \beta n$ and $\delta(G) \ge (1- {1 \over r}+\gamma)n$, then
$H \subset G$.
\end{conjecture}

The special case when $H$ is bipartite was shown by Abbasi~\cite{A}.
We will give an alternative proof of this by showing that if the band-width is small enough, then the graph is 
$\alpha$-separable
for a small enough $\alpha$. 

\begin{lemma}
Let $0<\beta <1$, and assume that $H$ is a graph of order $n$ with $bw(H) \le \beta n$. Then $H$ is a 
$\sqrt{\beta}$--separable graph.
\end{lemma}

\noindent {\bf Proof:}
We can decompose $H$ in the following way: Consider an ordering of the vertices of $H$ in which no
edge connects two vertices which are farther away from each other than $\beta n$. Divide the ordering
into $m=1 / \beta$ intervals. For simplicity we assume, that $m$ and $\sqrt{m}$ are integers and  $n$ 
is divisible by $m$. The $i$th interval, 
$I_i$ will contain the  vertices of order $(i-1)\beta n +1, \ldots,i \beta n$.

We let $$S=\bigcup_{i=1}^{\sqrt{m}}{I_{i\sqrt{m}}} \ ,$$ and for $0 \le j \le \sqrt{m}-1$, set 
$$C_j=\bigcup_{i=1}^{\sqrt{m}-1}{I_{j \sqrt{m} +i}} \ .$$ 

Clearly, $S$ consists of $\sqrt{m}$ intervals, each of length $n/m$, thus $|S| \le n/\sqrt{m}$. If $x \in C_j$ 
and $y \in C_k$ for
$j \neq k$, then $(x,y) \not\in E(H)$ because $bw(H) \le \beta n$. Hence, we have found a simple decomposition 
of $H$ which 
proves that $H$ is $\sqrt{\beta}$-separable.

\hfill \qedf

Unfortunately, if $\chi(H) \ge 3$, then our result does not imply Conjecture~\ref{BK}.

\medskip

%


\begin{thebibliography}{20}

\bibitem{A} S. Abbasi (1998),
Spanning subgraphs of dense graphs, Ph.D. theses, Department of Computer Science, Rutgers, the State University 
of New Jersey.

\bibitem{AY92} N. Alon and R. Yuster (1992),
Almost $H$-factors in dense graphs, Graphs and Combinatorics, {\bf 8}, 95--102.


\bibitem{AY96} N. Alon and R. Yuster (1996),
$H$-factors in dense graphs, Journal of Combinatorial Theory, Series B, {\bf 66}, 269--282.

\bibitem{B} B. Bollob\'as (1978), Extremal graph theory, Academic Press, London. 

\bibitem{Cs03} B. Csaba (2003), On the Bollob\'as--Eldridge conjecture for bipartite graphs,
submitted for publication.

\bibitem{ESi} P. Erd\H{o}s and M. Simonovits (1966), A limit theorem in graph theory, Studia Sci.
Math. Hungar., {\bf 1}, 51--57.

\bibitem{ESt} P. Erd\H{o}s and A. H. Stone (1946), On the structure of linear graphs, 
Bulletin of the  American Mathematical Society, {\bf 52}, 1089--1091.

\bibitem{HSz} A. Hajnal and E. Szemer\'edi (1970) Proof of a Conjecture
of Erd\H{o}s, in ``Combinatorial Theory and Its Applications, II''
(P. Erd\H{o}s, and V. T.  S\'os, Eds.), Colloquia Mathematica Societatis
J\'anos Bolyai, North-Holland, Amsterdam/London.

\bibitem{K} J. Koml\'os (1999), The Blow--up Lemma (survey), Combinatorics,
Probability and Computing, {\bf 8}, 161--176.

\bibitem{KSSz95} J. Koml\'os, G.N. S\'ark\"ozy and E. Szemer\'edi (1995), 
Proof of a packing conjecture of Bollob\'as,
Combinatorics, Probability and Computing, {\bf 4}, 241--255.


\bibitem{KSSz97} J. Koml\'os, G.N. S\'ark\"ozy and E. Szemer\'edi (1997) 
Blow-up Lemma, Combinatorica, {\bf 17}, 109-123.

\bibitem{KSSz98} J. Koml\'os, G.N. S\'ark\"ozy and E. Szemer\'edi (1998), An
Algorithmic Version of the Blow-up Lemma, Random Structures and Algorithms, {\bf 12}, 297-312.

\bibitem{KSSz} J. Koml\'os, G.N. S\'ark\"ozy and E. Szemer\'edi (2001), 
Proof of the Alon--Yuster conjecture, Discrete Mathematics, 255--269.

\bibitem{KSSz2} J. Koml\'os, G.N. S\'ark\"ozy and E. Szemer\'edi (2001), 
Spanning trees in dense graphs, Combinatorics, Probability and Computing, 397--416.


\bibitem{KS93} J. Koml\'os and M. Simonovits (1993), Szemer\'edi's Regularity Lemma and
its Applications in Graph Theory (survey), Combinatorics, Paul Erd\H{o}s is eighty,
Vol. 2 (Keszthely, 1993), 295--352.

\bibitem{Sz76} E. Szemer\'edi (1976), Regular Partitions of Graphs, Colloques
Internationaux C.N.R.S N$^{\underline {\rm o}}$ 260 - Probl\`emes
Combinatoires et Th\'eorie des Graphes, Orsay, 399-401.

\bibitem{T} P. Tur\'an (1941), On an extremal problem in graph theory (in Hungarian), 
Matematikai \'es Fizikai Lapok, {\bf 48}, 436--452.
\end{thebibliography}
\end{document}